# Poisson process approximation: From Palm theory to Stein's method

**Louis H. Y. Chen**[1,*] **and Aihua Xia**[2,†]

*National University of Singapore and University of Melbourne*

**Abstract:** This exposition explains the basic ideas of Stein's method for Poisson random variable approximation and Poisson process approximation from the point of view of the immigration-death process and Palm theory. The latter approach also enables us to define local dependence of point processes [Chen and Xia (2004)] and use it to study Poisson process approximation for locally dependent point processes and for dependent superposition of point processes.

## 1. Poisson approximation

Stein's method for Poisson approximation was developed by Chen [13] which is based on the following observation: a nonnegative integer valued random variable $W$ follows Poisson distribution with mean $\lambda$, denoted as $\text{Po}(\lambda)$, if and only if

$$\mathbb{E}\{\lambda f(W+1) - W f(W)\} = 0$$

for all bounded $f\colon \mathbb{Z}_+ \to \mathbb{R}$, where $\mathbb{Z}_+ := \{0, 1, 2, \ldots\}$. Heuristically, if $\mathbb{E}\{\lambda f(W+1) - W f(W)\} \approx 0$ for all bounded $f\colon \mathbb{Z}_+ \to \mathbb{R}$, then $\mathcal{L}(W) \approx \text{Po}(\lambda)$. On the other hand, as our interest is often on the difference $\mathbb{P}(W \in A) - \text{Po}(\lambda)(A) = \mathbb{E}[\mathbf{1}_A(W) - \text{Po}(\lambda)(A)]$, where $A \subset \mathbb{Z}_+$ and $\mathbf{1}_A$ is the indicator function on $A$, it is natural to relate the function $\lambda f(w+1) - w f(w)$ with $\mathbf{1}_A(w) - \text{Po}(\lambda)(A)$, leading to the Stein equation:

$$(1) \qquad \lambda f(w+1) - w f(w) = \mathbf{1}_A(w) - \text{Po}(\lambda)(A).$$

If the equation permits a bounded solution $f_A$, then

$$\mathbb{P}(W \in A) - \text{Po}(\lambda)(A) = \mathbb{E}\{\lambda f_A(W+1) - W f_A(W)\};$$

and

$$d_{TV}(\mathcal{L}(W), \text{Po}(\lambda)) := \sup_{A \subset \mathbb{Z}_+} |\mathbb{P}(W \in A) - \text{Po}(\lambda)(A)|$$
$$= \sup_{A \subset \mathbb{Z}_+} |\mathbb{E}\{\lambda f_A(W+1) - W f_A(W)\}|.$$

[1]Institute for Mathematical Sciences, National University of Singapore, 3 Prince George's Park, Singapore 118402, Republic of Singapore, e-mail: `matchyl@nus.edu.sg`
[2]Department of Mathematics and Statistics, University of Melbourne, Parkville, VIC 3010, Australia, e-mail: `xia@ms.unimelb.edu.au`
[*]Supported by research grant R-155-000-051-112 of the National University of Singapore.
[†]Supported by the ARC Centre of Excellence for Mathematics and Statistics of Complex Systems.
*AMS 2000 subject classifications:* primary 60–02; secondary 60F05, 60G55.
*Keywords and phrases:* immigration-death process, total variation distance, Wasserstein distance, locally dependent point process, hard core process, dependent superposition of point processes, renewal process.

As a special case in applications, we consider independent Bernoulli random variables $X_1, \cdots, X_n$ with $\mathbb{P}(X_i = 1) = 1 - \mathbb{P}(X_i = 0) = p_i$, $1 \leq i \leq n$, and $W = \sum_{i=1}^n X_i$, $\lambda = \mathbb{E}(W) = \sum_{i=1}^n p_i$. Since

$$\mathbb{E}[Wf(W)] = \sum_{i=1}^n \mathbb{E}[X_i f(W)] = \sum_{i=1}^n p_i \mathbb{E} f(W_i + 1),$$

where $W_i = W - X_i$, we have

$$\mathbb{E}\{\lambda f_A(W+1) - W f_A(W)\} = \sum_{i=1}^n p_i \mathbb{E}\left[f_A(W+1) - f_A(W_i + 1)\right]$$
$$= \sum_{i=1}^n p_i^2 \mathbb{E}\Delta f_A(W_i + 1),$$

where $\Delta f_A(i) = f_A(i+1) - f_A(i)$. Further analysis shows that $|\Delta f_A(w)| \leq \frac{1-e^{-\lambda}}{\lambda}$ (see [6] for an analytical proof and [26] for a probabilistic proof). Therefore

$$d_{TV}(\mathcal{L}(W), \text{Po}(\lambda)) \leq \left(1 \wedge \frac{1}{\lambda}\right) \sum_{i=1}^n p_i^2.$$

Barbour and Hall [7] proved that the lower bound of $d_{TV}(\mathcal{L}(W), \text{Po}(\lambda))$ above is of the same order as the upper bound. Thus this simple example of Poisson approximation demonstrates how powerful and effective Stein's method is. Furthermore, it is straightforward to use Stein's method to study the quality of Poisson approximation to the sum of dependent random variables which has many applications (see [18] or [8] for more information).

## 2. Poisson process approximation

Poisson process plays the central role in modeling the data on occurrence of rare events at random positions in time or space and is a building block for many other models such as Cox processes, marked Poisson processes (see [24]), compound Poisson processes and Lévy processes. To adapt the above idea of Poisson random variable approximation to Poisson process approximation, we need a probabilistic interpretation of Stein's method which was introduced by Barbour [4]. The idea is to split $f$ by defining $f(w) = g(w) - g(w-1)$ and rewrite the Stein equation (1) as

(2) $\quad \mathcal{A}g(w) := \lambda[g(w+1) - g(w)] + w[g(w-1) - g(w)] = \mathbf{1}_A(w) - \text{Po}(\lambda)(A),$

where $\mathcal{A}$ is the generator of an immigration-death process $Z_w(t)$ with immigration rate $\lambda$, unit per capita death rate, $Z_w(0) = w$, and stationary distribution $\text{Po}(\lambda)$. The solution to the Stein equation (2) is

(3) $\quad g_A(w) = -\int_0^\infty \mathbb{E}[\mathbf{1}_A(Z_w(t)) - \text{Po}(\lambda)(A)]dt.$

This probabilistic approach to Stein's method has made it possible to extend Stein's method to higher dimensions and process settings. To this end, let $\Gamma$ be a compact metric space which is the carrier space of the point processes being approximated. Suppose $d_0$ is a metric on $\Gamma$ which is bounded by 1 and $\rho_0$ is a pseudo-metric on $\Gamma$ which is also bounded by 1 but generates a weaker topology. We use $\delta_x$ to denote the point mass at $x$, let $\mathcal{X} := \{\sum_{i=1}^{k} \delta_{\alpha_i} : \alpha_1, \ldots, \alpha_k \in \Gamma, k \geq 1\}$, $\mathcal{B}(\mathcal{X})$ be the Borel $\sigma$–algebra generated by the weak topology ([23], pp. 168–170): a sequence $\{\xi_n\} \subset \mathcal{X}$ converges weakly to $\xi \in \mathcal{X}$ if $\int_\Gamma f(x)\xi_n(dx) \to \int_\Gamma f(x)\xi(dx)$ as $n \to \infty$ for all bounded continuous functions $f$ on $\Gamma$. Such topology can also be generated by the metric $d_1$ defined below (see [27], Proposition 4.2). A *point process* on $\Gamma$ is defined as a measurable mapping from a probability space $(\Omega, \mathcal{F}, \mathbb{P})$ to $(\mathcal{X}, \mathcal{B}(\mathcal{X}))$ (see [23], p. 13). We use $\Xi$ to stand for a point process on $\Gamma$ with finite intensity measure $\boldsymbol{\lambda}$ which has total mass $\lambda := \boldsymbol{\lambda}(\Gamma)$, where $\boldsymbol{\lambda}(A) = \mathbb{E}\Xi(A)$, for all Borel set $A \subset \Gamma$. Let $\text{Po}(\boldsymbol{\lambda})$ denote the distribution of a Poisson process on $\Gamma$ with intensity measure $\boldsymbol{\lambda}$.

Since a point process on $\Gamma$ is an $\mathcal{X}$-valued random element, the key step of extending Stein's method from one dimensional Poisson approximation to higher dimensions and process settings is, instead of considering $\mathbb{Z}_+$-valued immigration-death process, we now need an immigration-death process defined on $\mathcal{X}$. More precisely, by adapting (2), Barbour and Brown [5] define the Stein equation as

$$
\begin{aligned}
(4) \quad \mathcal{A}g(\xi) :&= \int_\Gamma [g(\xi + \delta_x) - g(\xi)]\boldsymbol{\lambda}(dx) + \int_\Gamma [g(\xi - \delta_x) - g(\xi)]\xi(dx) \\
&= h(\xi) - \text{Po}(\boldsymbol{\lambda})(h),
\end{aligned}
$$

where $\text{Po}(\boldsymbol{\lambda})(h) = \mathbb{E}h(\zeta)$ with $\zeta \sim \text{Po}(\boldsymbol{\lambda})$. The operator $\mathcal{A}$ is the generator of an $\mathcal{X}$-valued immigration-death process $Z_\xi(t)$ with immigration intensity $\boldsymbol{\lambda}$, unit per capita death rate, $Z_\xi(0) = \xi \in \mathcal{X}$, and stationary distribution $\text{Po}(\boldsymbol{\lambda})$. Its solution is

$$
(5) \quad g_h(\xi) = -\int_0^\infty \mathbb{E}[h(Z_\xi(t)) - \text{Po}(\boldsymbol{\lambda})(h)]dt,
$$

(see [5]).

To measure the error of approximation, we use Wasserstein pseudo-metric which has the advantage of allowing us to lift the carrier space to a bigger carrier space. Of course, other metrics such as the total variation distance can also be considered and the only difference is to change the set of test functions $h$. Let

$$
\rho_1\left(\sum_{i=1}^m \delta_{x_i}, \sum_{j=1}^n \delta_{y_j}\right) := \begin{cases} 1 & \text{if } m \neq n, \\ \min_\pi \frac{1}{m} \sum_{i=1}^m \rho_0(x_i, y_{\pi(i)}) & \text{if } m = n \geq 1, \\ 0 & \text{if } n = m = 0, \end{cases}
$$

where the minimum is taken over all permutations $\pi$ of $\{1, 2, \ldots, m\}$. Clearly, $\rho_1$ is a metric (resp. pseudo-metric) if $\rho_0$ is a metric (resp. pseudo-metric) on $\mathcal{X}$. Set

$$
\mathcal{H} = \{h \text{ on } \mathcal{X} : |h(\xi_1) - h(\xi_2)| \leq \rho_1(\xi_1, \xi_2) \text{ for all } \xi_1, \xi_2 \in \mathcal{X}\}.
$$

For point processes $\Xi_1$ and $\Xi_2$, define

$$
\rho_2(\mathcal{L}(\Xi_1), \mathcal{L}(\Xi_2)) := \sup_{h \in \mathcal{H}} |\mathbb{E}h(\Xi_1) - \mathbb{E}h(\Xi_2)|,
$$

then $\rho_2$ is a metric (resp. pseudo-metric) on the distributions of point processes if $\rho_1$ is a metric (resp. pseudo-metric). In summary, we defined a Wasserstein pseudo-metric on the distributions of point processes on $\Gamma$ through a pseudo-metric on $\Gamma$ as shown in the following chart:

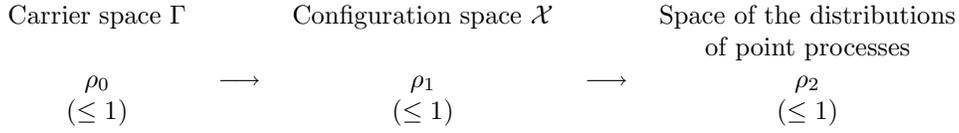

| Carrier space $\Gamma$ | Configuration space $\mathcal{X}$ | Space of the distributions of point processes |
|---|---|---|
| $\rho_0$ ($\leq 1$) | $\rho_1$ ($\leq 1$) | $\rho_2$ ($\leq 1$) |

As a simple example, we consider a Bernoulli process defined as

$$\Xi = \sum_{i=1}^{n} X_i \delta_{\frac{i}{n}},$$

where, as before, $X_1, \ldots, X_n$ are independent Bernoulli random variables with $\mathbb{P}(X_i = 1) = 1 - \mathbb{P}(X_i = 0) = p_i$, $1 \leq i \leq n$. Then $\Xi$ is a point process on carrier space $\Gamma = [0, 1]$ with intensity measure $\boldsymbol{\lambda} = \sum_{i=1}^{n} p_i \delta_{\frac{i}{n}}$. With the metric $\rho_0(x, y) = |x - y| := d_0(x, y)$, we denote the induced metric $\rho_2$ by $d_2$. Using the Stein equation (4), we have

$$\mathbb{E}h(\Xi) - \text{Po}(\boldsymbol{\lambda})(h)$$
$$= \mathbb{E}\left\{\int_\Gamma [g_h(\Xi + \delta_x) - g_h(\Xi)]\boldsymbol{\lambda}(dx) + \int_\Gamma [g_h(\Xi - \delta_x) - g_h(\Xi)]\Xi(dx)\right\}$$
$$= \sum_{i=1}^{n} p_i \mathbb{E}\left\{[g_h(\Xi + \delta_{\frac{i}{n}}) - g_h(\Xi)] - [g_h(\Xi_i + \delta_{\frac{i}{n}}) - g_h(\Xi_i)]\right\}$$
$$= \sum_{i=1}^{n} p_i^2 \mathbb{E}\left\{[g_h(\Xi_i + 2\delta_{\frac{i}{n}}) - g_h(\Xi_i + \delta_{\frac{i}{n}})] - [g_h(\Xi_i + \delta_{\frac{i}{n}}) - g_h(\Xi_i)]\right\},$$

where $\Xi_i = \Xi - X_i \delta_{\frac{i}{n}}$. It was shown in [27], Proposition 5.21, that

(6) $$\sup_{h \in \mathcal{H}, \alpha, \beta \in \Gamma} |g_h(\xi + \delta_\alpha + \delta_\beta) - g_h(\xi + \delta_\alpha) - g_h(\xi + \delta_\beta) + g_h(\xi)| \leq \frac{3.5}{\lambda} + \frac{2.5}{|\xi| + 1},$$

where, and in the sequel, $|\xi|$ is the total mass of $\xi$, $\lambda = \boldsymbol{\lambda}(\Gamma) = \sum_{i=1}^{n} p_i$. Hence

$$d_2(\mathcal{L}(\Xi), \text{Po}(\boldsymbol{\lambda})) = \sup_{h \in \mathcal{H}} |\mathbb{E}h(\Xi) - \text{Po}(\boldsymbol{\lambda})(h)|$$

(7) $$\leq \sum_{i=1}^{n} p_i^2 \left(\frac{3.5}{\lambda} + \mathbb{E}\frac{2.5}{\sum_{1 \leq j \leq n, j \neq i} X_j + 1}\right)$$
$$\leq \frac{6}{\lambda - \max_{1 \leq i \leq n} p_i} \sum_{i=1}^{n} p_i^2$$

since

$$\mathbb{E}\frac{1}{\sum_{1 \leq j \leq n, j \neq i} X_j + 1} = \mathbb{E}\int_0^1 z^{\sum_{1 \leq j \leq n, j \neq i} X_j} dz$$
$$= \int_0^1 \prod_{1 \leq j \leq n, j \neq i} [zp_j + (1 - p_j)] dz$$
$$\leq \int_0^1 \prod_{1 \leq j \leq n, j \neq i} e^{-p_j(1-z)} dz = \int_0^1 e^{-(\lambda - p_j)(1-z)} dz \leq \frac{1}{\lambda - p_j},$$

(see [27], pp. 167–168). Since $d_2(\mathcal{L}(\Xi), \text{Po}(\boldsymbol{\lambda})) \geq d_{TV}(\mathcal{L}(|\Xi|), \text{Po}(\lambda))$ and the lower bound of $d_{TV}(\mathcal{L}(|\Xi|), \text{Po}(\lambda))$ is of the same order as $\frac{1}{\lambda} \sum_{i=1}^{n} p_i^2$ [7], the bound in (7) is of the optimal order.

## 3. From Palm theory to Stein's method

Barbour's probabilistic approach to Stein's method is based on the conversion of a first order difference equation to a second order difference equation. In this section, we take another approach to Stein's method from the point of Palm theory. The connection between Stein's method and Palm theory has been known to many others (e.g., T. C. Brown (personnel communication), [9]) and the exposition here is mainly based on [14] and [27].

There are two properties which distinguish a Poisson process from other processes: independent increments and the number of points on any bounded set follows Poisson distribution. Hence, a Poisson process can be thought as a process pieced together by lots of independent "Poisson components" (if the location is an atom, the "component" will be a Poisson random variable, but if the location is diffuse, then the "component" is either 0 or 1) ([27], p. 121). Consequently, to specify a Poisson process $N$, it is sufficient to check that "each component" $N(d\alpha)$ is Poisson and independent of the others, that is $\mathbb{E}\{[\mathbb{E}N(d\alpha)]g(N+\delta_\alpha) - N(d\alpha)g(N)\} = 0$, which is equivalent to

$$(8) \qquad \frac{\mathbb{E}[g(N)N(d\alpha)]}{\mathbb{E}N(d\alpha)} = \mathbb{E}g(N+\delta_\alpha),$$

for all bounded function $g$ on $\mathcal{X}$ and all $\alpha \in \Gamma$ (see [27], p. 121). To make the heuristic argument rigorous, one needs the tools of Campbell measures and Radon-Nikodym derivatives ([23], p. 83).

In general, for each point process $\Xi$ with finite mean measure $\boldsymbol{\lambda}$, we may define the Campbell measure $C(B, M) = \mathbb{E}[\Xi(B)\mathbf{1}_{\Xi \in M}]$ for all Borel $B \subset \Gamma$, $M \in \mathcal{B}(\mathcal{X})$. This measure is finite and admits the following disintegration:

$$(9) \qquad C(B, M) = \int_B Q_s(M)\boldsymbol{\lambda}(ds),$$

or equivalently,

$$Q_s(M) = \frac{\mathbb{E}[\Xi(ds)\mathbf{1}_{\Xi \in M}]}{\boldsymbol{\lambda}(ds)}, \ M \in \mathcal{B}(\mathcal{X}), \ s \in \Gamma \ \boldsymbol{\lambda} \ a.s.,$$

where $\{Q_s, \ s \in \Gamma\}$ are probability measures on $\mathcal{B}(\mathcal{X})$ ([23], p. 83 and p. 164) and are called *Palm distributions*. Moreover, (9) is equivalent to that, for any measurable function $f: \Gamma \times \mathcal{X} \to \mathbb{R}_+$,

$$(10) \qquad \mathbb{E}\left(\int_B f(\alpha, \Xi)\Xi(d\alpha)\right) = \int_B \int_\mathcal{X} f(\alpha, \xi)Q_\alpha(d\xi)\boldsymbol{\lambda}(d\alpha)$$

for all Borel set $B \subset \Gamma$. A point process $\Xi_\alpha$ (resp. $\Xi_\alpha - \delta_\alpha$) on $\Gamma$ is called a *Palm process* (resp. *reduced Palm process*) of $\Xi$ at location $\alpha$ if it has the Palm distribution $Q_\alpha$ and, when $\Xi$ is a simple point process (a point process taking values 0 or 1 at each location), the Palm distribution $\mathcal{L}(\Xi_\alpha)$ can be interpreted as the conditional distribution of $\Xi$ given that there is a point of $\Xi$ at $\alpha$. It follows from (10) that the Palm process satisfies

$$\mathbb{E}\int_\Gamma f(\alpha, \Xi)\Xi(d\alpha) = \mathbb{E}\int_\Gamma f(\alpha, \Xi_\alpha)\boldsymbol{\lambda}(d\alpha)$$

for all bounded measurable functions $f$ on $\Gamma \times \mathcal{X}$. In particular, $\Xi$ is a Poisson process if and only if

$$\mathcal{L}(\Xi_\alpha) = \mathcal{L}(\Xi + \delta_\alpha), \ \boldsymbol{\lambda} \ a.s.$$

where the extra point $\delta_\alpha$ is due to the "Poisson property" of $\Xi\{\alpha\}$, and $\Xi_\alpha|_{\Gamma\setminus\{\alpha\}}$ has the same distribution as $\Xi|_{\Gamma\setminus\{\alpha\}}$ because of independent increments. Here $\xi|_A$ stands for the point measure restricted to $A \subset \Gamma$ ([23], p. 12). In other words, $\Xi \sim \text{Po}(\boldsymbol{\lambda})$ if and only if

$$\mathbb{E}\left\{\int_\Gamma f(\alpha, \Xi + \delta_\alpha)\boldsymbol{\lambda}(d\alpha) - \int_\Gamma f(\alpha, \Xi)\Xi(d\alpha)\right\} = 0,$$

for a sufficiently rich class of functions $f$, so we define

$$Df(\xi) := \int_\Gamma f(x, \xi + \delta_x)\boldsymbol{\lambda}(dx) - \int_\Gamma f(x, \xi)\xi(dx).$$

If $\mathbb{E}Df(\Xi) \approx 0$ for an appropriate class of test functions $f$, then $\mathcal{L}(\Xi_\alpha)$ is close to $\mathcal{L}(\Xi + \delta_\alpha)$, which means that $\mathcal{L}(\Xi)$ is close to $\text{Po}(\boldsymbol{\lambda})$ under the metric or pseudometric specified by the class of test functions $f$.

If $f_g$ is a solution of

$$Df(\xi) = g(\xi) - \text{Po}(\boldsymbol{\lambda})(g),$$

then a distance between $\mathcal{L}(\Xi)$ and $\text{Po}(\boldsymbol{\lambda})$ is measured by $|\mathbb{E}Df_g(\Xi)|$ over the class of functions $g$.

From above analysis, we can see that there are many possible solutions $f_g$ for a given function $g$. The one which admits an immigration-death process interpretation is by setting

$$f(x, \xi) = h(\xi) - h(\xi - \delta_x),$$

so that $Df$ takes the following form:

$$Df(\xi) = \int_\Gamma [h(\xi + \delta_x) - h(\xi)]\boldsymbol{\lambda}(dx) + \int_\Gamma [h(\xi - \delta_x) - h(\xi)]\xi(dx) = \mathcal{A}h(\xi),$$

where $\mathcal{A}$ is the same as the generator defined in section 2.

## 4. Locally dependent point processes

We say a point process $\Xi$ is *locally dependent* with neighborhoods $\{A_\alpha \subset \Gamma: \ \alpha \in \Gamma\}$ if $\mathcal{L}(\Xi|_{A_\alpha^c}) = \mathcal{L}(\Xi_\alpha|_{A_\alpha^c})$, $\alpha \in \Gamma$ $\boldsymbol{\lambda}$ a.s.

The following theorem is virtually from Corollary 3.6 in [14] combined with the new estimates of Stein's factors in [27], Proposition 5.21.

**Theorem 4.1.** *If $\Xi$ is a point process on $\Gamma$ with finite intensity measure $\boldsymbol{\lambda}$ which has the total mass $\lambda$ and locally dependent with neighborhoods $\{A_\alpha \subset \Gamma: \ \alpha \in \Gamma\}$. Then*

$$\rho_2(\mathcal{L}(\Xi), \text{Po}(\boldsymbol{\lambda})) \leq \mathbb{E}\int_{\alpha\in\Gamma} \left(\frac{3.5}{\lambda} + \frac{2.5}{|\Xi^{(\alpha)}| + 1}\right)(\Xi(A_\alpha) - 1)\Xi(d\alpha)$$

$$+ \mathbb{E}\int_{\alpha\in\Gamma}\int_{\beta\in A_\alpha}\left(\frac{3.5}{\lambda} + \frac{2.5}{|\Xi^{(\alpha)}_\beta| + 1}\right)\boldsymbol{\lambda}(d\alpha)\boldsymbol{\lambda}(d\beta),$$

*where $\Xi^{(\alpha)} = \Xi|_{A_\alpha^c}$ and $\Xi^{(\alpha)}_\beta = \Xi_\beta|_{A_\alpha^c}$.*

**Remark.** The error bound is a "correct" generalization of $\frac{1}{\lambda}\sum_{i=1}^n p_i^2$ with the Stein factor $\frac{1}{\lambda}$ replaced by a nonuniform bound.

## 5. Applications

### 5.1. Matérn hard core process on $\mathbb{R}^d$

A Matérn hard core process $\Xi$ on compact $\Gamma \subset \mathbb{R}^d$ is a model for particles with repulsive interaction. It assumes that points occur according to a Poisson process with uniform intensity measure on $\Gamma$. The configurations of $\Xi$ are then obtained by deleting any point which is within distance $r$ of another point, irrespective of whether the latter point has itself already been deleted [see Cox & Isham [17], p. 170].

The point process is locally dependent with neighborhoods $\{B(\alpha, 2r)\colon \alpha \in \Gamma\}$, where $B(\alpha, s)$ is the ball centered at $\alpha$ with radius $s$. Let $\boldsymbol{\lambda}$ be the intensity measure of $\Xi$, $d_0(\alpha, \beta) = \min\{|\alpha - \beta|, 1\}$, then

$$d_2(\mathcal{L}(\Xi), \mathrm{Po}(\boldsymbol{\lambda})) = O\left(\frac{\mu \mathrm{Vol}(B(0,1))(2r)^d}{\mathrm{Vol}(\Gamma)}\right),$$

where $\mu$ is the mean of the total number of points of the original Poisson process (see [14], Theorem 5.1).

### 5.2. Palindromes in a genome

Let $\{I_i\colon 1 \le i \le n\}$ be locally dependent Bernoulli random variables, $\{U_i\colon 1 \le i \le n\}$ be independent $\Gamma$-valued random elements which are also independent of $\{I_i\colon 1 \le i \le n\}$, set $\Xi = \sum_{i=1}^n I_i \delta_{U_i}$, then $\Xi$ is a point process on $\Gamma$. For $U_i = i/n$ this point process models palindromes in a genome where $I_i$ represents whether a palindrome occurs at $i/n$. The point process can also be used to describe the vertices in a random graph.

In general, the $U_i$'s could take the same value and one cannot tell which $U_i$ and therefore which $I_i$ contributes to the value. To overcome this difficulty we lift the process up to a point process $\Xi' = \sum_{i=1}^n I_i \delta_{(i, U_i)}$ on a larger space $\Gamma' = \{1, 2, \ldots, n\} \times \Gamma$. The metric $d_0$ becomes a pseudo-metric $\rho_0$, that is, $\rho_0((i, s), (j, t)) = d_0(s, t)$, and $\Xi'$ a locally dependent process (see [14], section 4). It turns out that the Poisson process approximation of $\Xi = \sum_{i=1}^n I_i \delta_{U_i}$ is a special case of the following section.

### 5.3. Locally dependent superposition of point processes

Since the publication of the Grigelionis Theorem [20] which states that the superposition of independent sparse point processes on carrier space $\mathbb{R}_+$ is close to a Poisson process, there has been a lot of study on the weak convergence of point processes to a Poisson process under various conditions (see, e.g., [16, 19, 21] and [10]). Extensions to dependent superposition[1] of sparse point processes have been carried out in [1, 2, 3, 11, 22]. Schuhmacher [25] considered the Wasserstein distance between the weakly dependent superposition of sparse point processes and a Poisson process.

Let $\Gamma$ be a compact metric space, $\{\Xi_i\colon i \in \mathcal{I}\}$ be a collection of point processes on $\Gamma$ with intensity measures $\boldsymbol{\lambda}_i$, $i \in \mathcal{I}$. Define $\Xi = \sum_{i \in \mathcal{I}} \Xi_i$ with intensity measure $\boldsymbol{\lambda} = \sum_{i \in \mathcal{I}} \boldsymbol{\lambda}_i$. Assume $\{\Xi_i \colon i \in \mathcal{I}\}$ are locally dependent: that is, for each $i \in \mathcal{I}$, there exists a neighbourhood $A_i \subset \mathcal{I}$ such that $i \in A_i$ and $\Xi_i$ is independent of $\{\Xi_j \colon j \notin A_i\}$.

---

[1]We use "(resp. locally, weakly) dependent superposition of point processes" to mean that the point processes are (resp. locally, weakly) dependent among themselves.

The locally dependent point process $\Xi = \sum_{i=1}^n I_i \delta_{U_i}$ can be regarded as a locally dependent superposition of point processes defined above.

**Theorem 5.1** ([15]). *With the above setup, $\lambda = \boldsymbol{\lambda}(\Gamma)$, we have*

$$d_2(\mathcal{L}(\Xi), \mathrm{Po}(\boldsymbol{\lambda})) \leq \mathbb{E} \sum_{i \in \mathcal{I}} \left( \frac{3.5}{\lambda} + \frac{2.5}{|\Xi^{(i)}| + 1} \right) \int_\Gamma d_1'(\mathbf{V}_i, \mathbf{V}_{i,\alpha}) \boldsymbol{\lambda}_i(d\alpha)$$
$$+ \sum_{i \in \mathcal{I}} \left( \frac{3.5}{\lambda} + \mathbb{E} \frac{2.5}{|\Xi^{(i)}| + 1} \right) \mathbb{E} \int_\Gamma d_1'(\Xi_i, \Xi_{i,\alpha}) \boldsymbol{\lambda}_i(d\alpha),$$

*where $\Xi^{(i)} = \sum_{j \notin A_i} \Xi_j$, $\mathbf{V}_i = \sum_{j \in A_i \setminus \{i\}} \Xi_j$, $\Xi_{i,\alpha}$ is the reduced Palm process of $\Xi_i$ at $\alpha$,*

$$\mathbb{P}(\mathbf{V}_{i,\alpha} \in M) = \frac{\mathbb{E}[\Xi_i(d\alpha) \mathbf{1}_{\mathbf{V}_i \in M}]}{\mathbb{E} \Xi_i(d\alpha)} \text{ for all } M \in \mathcal{B}(\mathcal{X})$$

*and*

$$d_1'(\xi_1, \xi_2) = \min_{\pi \colon \text{ permutations of } \{1,\ldots,m\}} \sum_{i=1}^n d_0(y_i, z_{\pi(i)}) + (m - n)$$

*for $\xi_1 = \sum_{i=1}^n \delta_{y_i}$ and $\xi_2 = \sum_{i=1}^m \delta_{z_i}$ with $m \geq n$ [Brown & Xia [12]].*

**Corollary 5.2** ([14]). *For $\Xi = \sum_{i \in \mathcal{I}} I_i \delta_{U_i}$ and $\lambda = \sum_{i \in \mathcal{I}} p_i$ defined in section 5.2,*

$$d_2(\mathcal{L}(\Xi), \mathrm{Po}(\boldsymbol{\lambda})) \leq \mathbb{E} \sum_{i \in \mathcal{I}} \sum_{j \in A_i \setminus \{i\}} \left( \frac{3.5}{\lambda} + \frac{2.5}{V_i + 1} \right) I_i I_j$$
$$+ \sum_{i \in \mathcal{I}} \sum_{j \in A_i} \left( \frac{3.5}{\lambda} + \mathbb{E} \left[ \frac{2.5}{V_i + 1} \middle| I_j = 1 \right] \right) p_i p_j,$$

*where $V_i = \sum_{j \notin A_i} I_j$.*

**Corollary 5.3** ([15]). *Suppose that $\{\Xi_i \colon 1 \leq i \leq n\}$ are independent renewal processes on $[0, T]$ with the first arrival time of $\Xi_i$ having distribution $G_i$ and its inter-arrival time having distribution $F_i$, and let $\Xi = \sum_{i \in \mathcal{I}} \Xi_i$ and $\boldsymbol{\lambda}$ be its intensity measure, then*

$$d_2(\mathcal{L}(\Xi), \mathrm{Po}(\boldsymbol{\lambda})) \leq \frac{6 \sum_{i=1}^n [2F_i(T) + G_i(T)] G_i(T)/(1 - F_i(T))^2}{\sum_{i=1}^n G_i(T) - \max_j \frac{G_j(T)}{1 - F_j(T)}}.$$

**References**

[1] BANYS, R. (1975). The convergence of sums of dependent point processes to Poisson processes. (Russian. Lithuanian, English summary) *Litovsk. Mat. Sb.* **15** 11–23, 223. MR0418231

[2] BANYS, R. (1985). A Poisson limit theorem for rare events of a discrete random field. (Russian. English, Lithuanian summary) *Litovsk. Mat. Sb.* **25** 3–8. MR0795851

[3] BANYS, R. (1980). On superpositions of random measures and point processes. *Mathematical Statistics and Probability Theory.* Proc. Sixth Internat. Conf., Wisla, 1978, *Lecture Notes in Statist.* **2**. Springer, New York-Berlin, pp. 26–37. MR0577268

[4] BARBOUR, A. D. (1988). Stein's method and Poisson process convergence. *J. Appl. Probab.* **25** (A) 175–184. MR0974580

[5] BARBOUR, A. D. AND BROWN, T. C. (1992). Stein's method and point process approximation. *Stoch. Procs. Applics* **43** 9–31. MR1190904

[6] BARBOUR, A. D. AND EAGLESON, G. K. (1983). Poisson approximation for some statistics based on exchangeable trials. *Adv. Appl. Prob.* **15** 585–600. MR0706618

[7] BARBOUR, A. D. AND HALL, P. (1984). On the rate of Poisson convergence. *Math. Proc. Cambridge Philos. Soc.* **95** 473–480. MR0755837

[8] BARBOUR, A. D., HOLST, L. AND JANSON, S. (1992). *Poisson Approximation.* Oxford Univ. Press. MR1163825

[9] BARBOUR, A. D. AND MÅNSSON, M. (2002). Compound Poisson process approximation. *Ann. Probab.* **30** 1492–1537. MR1920275

[10] BROWN, T. C. (1978). A martingale approach to the Poisson convergence of simple point processes. *Ann. Probab.* **6** 615–628. MR0482991

[11] BROWN, T. C. (1979). Position dependent and stochastic thinning of point processes. *Stoch. Procs. Applics* **9** 189–193. MR0548838

[12] BROWN, T. C. AND XIA, A. (1995). On metrics in point process approximation. *Stochastics and Stochastics Reports* **52** 247–263. MR1381671

[13] CHEN, L. H. Y. (1975). Poisson approximation for dependent trials. *Ann. Probab.* **3** 534–545. MR0428387

[14] CHEN, L. H. Y. AND XIA, A. (2004). Stein's method, Palm theory and Poisson process approximation. *Ann. Probab.* **32** 2545–2569. MR2078550

[15] CHEN, L. H. Y. AND XIA, A. (2006). Poisson process approximation for dependent superposition of point processes. *(preprint).*

[16] ÇINLAR, E. (1972). Superposition of point processes. *Stochastic Point Processes: Statistical Analysis, Theory, and Applications.* Conf., IBM Res. Center, Yorktown Heights, NY, 1971. Wiley-Interscience, New York, pp. 549–606. MR0365697

[17] COX, D. R. AND ISHAM, V. (1980). *Point Processes.* Chapman & Hall. MR0598033

[18] ERHARDSSON, T. (2005). Poisson and compound Poisson approximation. In: *An Introduction to Stein's Method*, Eds. A. D. Barbour and L. H. Y. Chen. World Scientific Press, Singapore, pp. 61–113. MR2235449

[19] GOLDMAN, J. R. (1967). Stochastic point processes: limit theorems. *Ann. Math. Statist.* **38** 771–779. MR0217854

[20] GRIGELIONIS, B. (1963). On the convergence of sums of random step processes to a Poisson process. *Theor. Probab. Appl.* **8** 177–182. MR0152013

[21] JAGERS, P. (1972). On the weak convergence of superpositions of point processes. *Z. Wahrsch. verw. Geb.* **22** 1–7. MR0309191

[22] KALLENBERG, O. (1975). Limits of compound and thinned point processes. *J. Appl. Probab.* **12** 269–278. MR0391251

[23] KALLENBERG, O. (1983). *Random Measures.* Academic Press, London. MR0818219

[24] KINGMAN, J. F. C. (1993). *Poisson processes.* Oxford University Press. MR1207584

[25] SCHUHMACHER, D. (2005). Distance estimates for dependent superpositions of point processes. *Stoch. Procs. Applics* **115** 1819–1837. MR2172888

[26] XIA, A. (1999). A probabilistic proof of Stein's factors. *J. Appl. Probab.* **36** 287–290. MR1699603

[27] XIA, A. (2005). Stein's method and Poisson process approximation. In: *An*

*Introduction to Stein's Method*, Eds. A. D. Barbour and L. H. Y. Chen. World Scientific Press, Singapore, pp. 115–181. MR2235450